\documentclass[12pt]{article}

\usepackage{amsmath}
\usepackage{amsfonts}
\usepackage{latexsym}
\usepackage{amssymb}
\usepackage{curves}
\usepackage{graphicx}
\usepackage{amsthm,amscd, amssymb}

\newtheorem{thm}{Theorem}[section]
\newtheorem{pro}[thm]{Proposition}
\newtheorem{cor}[thm]{Corollary}
\newtheorem{lem}[thm]{Lemma}

\newtheorem{rem}[thm]{Remark}

\newcommand{\mapsfrom}
{\mathrel{\reflectbox{\ensuremath{\mapsto}}}}

\newcommand{\back}{\backslash}

\newcommand{\ua}{\uparrow}

\newcommand{\smid}{\, | \,}

\newcommand{\noin}{\noindent}

\newcommand{\tinyskip}{\vspace{0.05in}}

\newcommand{\FF}{{\mathbb{F}}}

\newcommand{\PP}{{\mathbb{P}}}

\renewcommand{\SS}{{\mathbb{S}}}

\newcommand{\Dia}{\mbox{\rm{Dia}}}
\newcommand{\End}{\mbox{\rm{End}}}

\newcommand{\Ind}{\mbox{\rm{Ind}}}
\newcommand{\Inf}{\mbox{\rm{Inf}}}

\newcommand{\Irr}{\mbox{\rm{Irr}}}
\newcommand{\Iso}{\mbox{\rm{Iso}}}

\newcommand{\Res}{\mbox{\rm{Res}}}

\newcommand{\br}{\mbox{\rm{br}}}

\newcommand{\id}{\mbox{\rm{id}}}

\newcommand{\sspan}{\mbox{\rm{span}}}
\newcommand{\tr}{\mbox{\rm{tr}}}

\newcommand{\oo}{\overline}

\newcommand{\ooN}{{\overline{N}}}

\newcommand{\cL}{{\cal L}}

\newcommand{\cP}{{\cal P}}

\begin{document}


\title{The pointed $p$-groups on a block algebra}

\author{\large Laurence Barker \\ \mbox{} \\
\normalsize Department of Mathematics \\
\normalsize Bilkent University \\
\normalsize 06800 Bilkent, Ankara \\
\normalsize Turkey}

\maketitle

\small

\begin{abstract}
\noin A pointed $p$-group is a pointed group
$P_\gamma$ such that $P$ is a $p$-group.
We parameterize the pointed $p$-groups on
a group algebra or on a block algebra of a
group algebra. The parameterization involves
$p$-subgroups and irreducible characters of
centralizers of $p$-subgroups.

\smallskip
\noin 2020 {\it Mathematics Subject Classification:} 20C20.

\smallskip
\noin {\it Keywords:} Puig category; relative multiplicity;
pointed $p$-group; generalized piece.
\end{abstract}

\section{Introduction}
\label{1}

For a block of a group algebra, Puig \cite{Pui86}
introduced a refinement of the fusion system
which Th\'{e}venaz \cite[Section 47]{The95}
called the Puig category of the block. Let
$\FF$ be an algebraically closed field of
prime characteristic $p$, let $G$ be a
finite group and let $b$ be a block of the
group algebra $\FF G$. Let $D_\lambda$
be a maximal local pointed group on the
block algebra $\FF G b$. Recall that the
{\bf Puig category} $\cL$ of $\FF G b$
associated with $D_\lambda$ is defined as
follows. The objects of $\cL$ are the
local pointed subgroups of $D_\lambda$.
For local pointed subgroups $P_\gamma$
and $Q_\delta$ of $D_\lambda$, the
$\cL$-morphisms $P_\gamma \leftarrow
Q_\delta$ are the conjugation monomorphisms
$P \ni {}^g y \mapsfrom y \in Q$ where
$g \in G$ satisfies $P_\gamma \geq
{}^g (Q_\delta)$. The composition is the
usual composition of group monomorphisms.
Since $D_\lambda$ is determined by
$\FF G b$ up to $G$-conjugacy, the
category $\cL$ is determined by $\FF G b$
up to isomorphism of categories.

To explicitly specify a local pointed group
$P_\gamma$ on $\FF G b$, or more generally
on $\FF G$, we make use of Brauer characters.
We understand an $\FF G$-character to be a
Brauer character of $\FF G$, defined by means
of a fixed embedding of the torsion unit group
of $\FF$ in the torsion unit group of some
algebraically closed field of characteristic $0$.
A well-known result of Puig, recorded as Theorem
\ref{2.1} below, expresses, for any $p$-subgroup
$P$ of $G$, a bijective correspondence between
the local points of $P$ on $\FF G$ and the
irreducible $\FF C_G(P)$-characters. However,
the Puig category is hard to determine explicitly
because the inclusion relation between the local
pointed groups on $\FF G$ is difficult to describe
in terms of the corresponding irreducible
characters.

We shall assume familiarity with the theory
of $G$-algebras. For accounts of the topic, see
Linckelmann \cite{Lin18}, Th\'{e}venaz
\cite{The95}. For any pointed groups $W_\omega$
and $U_\mu$ on $\FF G$ with $W \leq U$,
write $m(W_\omega, U_\mu)$ to denote the
relative multiplicity of $W_\omega$ in
$U_\mu$. Recall that $W_\omega \leq U_\mu$
if and only if $m(W_\omega, U_\mu) \neq 0$.
It is easy to see that, given $W \leq V \leq U$,
then we have a matrix relation
$$m(W_\omega, U_\mu) = \sum_\nu
  m(W_\omega, V_\nu) \, m(V_\nu, U_\mu)$$
where $\nu$ runs over the points of $V$ on
$\FF G$. That matrix relation cannot, in general,
be confined to local pointed groups. Indeed,
supposing that $W_\omega$ and $U_\mu$ are
local pointed groups on $\FF G$, whereupon
$W$, $V$, $U$ must be $p$-groups, evaluation
of the sum still requires us to consider all the
points $\nu$ of $V$ on $\FF G$, not just the
local points of $V$ on $\FF G$. 
We call a pointed group $U_\mu$ on $\FF G$
a {\bf pointed $p$-group} when
$U$ is a $p$-group. In order to make use
of the above matrix relation, as was done
in the proof of [BG22, 5.2] for instance, it may
be necessary to consider all the pointed
$p$-groups on $\FF G$ or on $\FF G b$.

The main aim of this paper is to parameterize
the pointed $p$-groups on $\FF G$ and,
in particular, on $\FF G b$. That is to say,
we shall put the pointed $p$-groups in a
bijective correspondence with the elements
of a set that can be explicitly determined.

In Section \ref{2}, we shall recall a result
of Puig which establishes a bijective
correspondence between the local pointed
groups on $\FF G$ and some pairs which
we call the pieces of $\FF G$. We shall
define the notion of a piece in Section
\ref{2}. For now, let us just say that, if
one is armed with an explicit description
of the poset of $p$-subgroups of $G$ and
the modular character tables of the
centralizers of the $p$-subgroups of $G$,
then one knows explicitly what the pieces of
$\FF G$ are. Towards a parameterization of
the pointed $p$-subgroups on $\FF G$,
we shall introduce the notion of a generalized
piece of $\FF G$. Again, granted an explicit
description of the poset of $p$-subgroups
and the modular character tables, one knows
explicitly what the generalized pieces of
$\FF G$ are. Each pointed $p$-group
on $\FF G$ is associated with a unique
generalized piece of $\FF G$. The question
to be answered, then, is as to which of the
generalized pieces are associated with
local pointed groups. We shall introduce
the notion of a substantive generalized piece.
The defining condition for substantivity
involves involves a simple module
constructed for a semidirect product
$C_G(Q) {\rtimes} (N_P(Q)/Q)$, where
$Q \leq P$ are $p$-subgroups of $G$.
Proposition \ref{2.3} allows the condition
to be reformulated in terms of a
Clifford-theoretic construction.

Our main result, Theorem \ref{2.8}, stated
in Section \ref{2}, describes a bijective
correspondence between the pointed
$p$-groups on $\FF G$ and the substantive
generalized pieces of $\FF G$. Corollary
\ref{2.9} describes how the bijective
correspondence is compatible with blocks.

The proof of Theorem \ref{2.8} will be given
in Section 3, where we shall introduce
the notion of the absolute multiplicity
of a generalised piece. The substantive
generalized pieces are those generalized
pieces whose absolute multiplicity is
nonzero. Proposition \ref{3.12} says
that the bijective correspondence
between the local pointed groups
and the substantive generalized
pieces preserves absolute multiplicity.

To illustrate the theory, we shall present
two examples in Section 4. For the
principal $2$-blocks of the symmetric
groups $S_4$ and $S_5$, we shall
calculate the relative multiplicities between
the substantive generalized pieces. That
will yield, in particular, the relative
multiplicities between the pieces of the
blocks, in other words, the relative
multiplicities between the local
pointed groups. The method is to first
determine which generalized pieces
are substantive (in effect, classifying
the pointed $p$-groups on the block
algebras), then determining the
absolute multiplicities of the
substantive generalized pieces
(in effect, determining the absolute
multiplicities of the pointed $p$-groups).

\section{Qualitative results}
\label{2}

In this section, we shall define the
notions of a piece of $\FF G$, a
generalized piece of $\FF G$ and a
substantive generalized piece of
$\FF G$. In Theorem \ref{2.8}, we shall
describe a bijective correspondence
between the pointed $p$-groups on
$\FF G$ and the substantive generalized
pieces of $\FF G$. Three of the results
in this section, Propositions \ref{2.3},
\ref{2.7} and Theorem \ref{2.8}, are
expressed in a qualitative way, and
their proofs, which rely on some
formulas for multiplicities, will be
deferred to the next section.

We let $\Irr(\FF G)$ denote the set of
irreducible $\FF G$-characters. Given
$\xi \in \Irr(\FF G)$, we write $V(\xi)$
to denote a simple $\FF G$-module with
modular character $\xi$. We write $E(\xi)$
to denote an indecomposable projective
$\FF G$-module with a quotient isomorphic
to $V(\xi)$. Of course, $V(\xi)$ and $E(\xi)$
are well-defined up to isomorphism.

For a $p$-subgroup $P$ of $G$, we write
the $P$-relative Brauer map as
$\br_P : \FF C_G(P) \leftarrow (\FF G)^P$.
We define a {\bf piece} of $\FF G$ to be
a pair $(P, \theta)$, usually written as
$P_\theta$, where $P$ is a $p$-subgroup
of $G$ and $\theta \in \Irr(\FF C_G(P))$.
We allow $G$ to permute the pieces of
$\FF G$ by defining ${}^g (P_\theta) =
({}^g P)_{{}^g \theta}$ for $g \in G$. In
Theorem \ref{2.8} below, we shall be extending
the following fundamental theorem, which can
be found in Th\'{e}venaz \cite[37.6]{The95}.

\begin{thm} \label{2.1}
{\rm (Puig.)} There is a $G$-equivariant
bijective correspondence between:

\tinyskip
\noin $\bullet$ the local pointed groups
$P_\gamma$ on $\FF G$,

\tinyskip
\noin $\bullet$ the pieces $P'_\theta$ of
$\FF G$,

\tinyskip
\noin such that $P_\gamma \leftrightarrow
P'_\theta$ if and only if $P = P'$ and
$E(\theta) \cong \FF C_G(P) \br_P(i)$
where $i \in \gamma$. The condition is
independent of the choice of $i$.
\end{thm}

To generalize that theorem, we shall
need to generalize the notion of a piece.
Consider the pairs $(P, Q_\phi)$
where $P$ is a $p$-subgroup of $G$
and $Q_\phi$ is a piece of $\FF G$
such that $P \geq Q$. Two such pairs
$(P, Q_\phi)$ and $(P', Q'_{\phi'})$ are
to be deemed equivalent provided
$P = P'$ and the pieces $Q_\phi$ and
$Q'_{\phi'}$ are $P$-conjugate. We write
$P {\ua} Q_\phi$ to denote the equivalence
class of $(P, Q_\phi)$. We call
$P {\ua} Q_\phi$ a {\bf generalized piece}
of $\FF G$. We allow $G$ to permute the
generalized pieces of $\FF G$ by defining
${}^g(P {\ua} Q_\phi) = ({}^g P) {\ua}
{}^g (Q_\phi)$ for $g \in G$.
We identify any piece $P_\theta$ with the
generalized piece $P {\ua} P_\theta$.

To define the substantivity condition on
generalized pieces, we shall be needing
the following abstract lemma.

\begin{lem} \label{2.2}
Let $K \unlhd G$ such that $G/K$ is a
$p$-group. Given $\xi \in \Irr(\FF G)$
and $\eta \in \Irr(\FF K)$, then
$E(\xi) \cong {}_G \Ind {}_K (E(\eta))$
if and only if ${}_K \Res {}_G (V(\xi))$
is a direct sum of mutually
non-isomorphic $G$-conjugates of
$V(\eta)$. Furthermore, those
equivalent conditions characterize
a bijective correspondence
$\xi \leftrightarrow [\eta]_G$ between
the irreducible $\FF G$-characters
$\xi$ and the $G$-conjugacy classes
$[\eta]_G$ of irreducible
$\FF K$-characters $\eta$.
\end{lem}

\begin{proof}
A special case of Linckelmann
\cite[5.12.10]{Lin18} asserts that any
primitive idempotent of $\FF K$
remains primitive in $\FF G$. So
there is a function $\Irr(\FF K) \ni
\eta \mapsto \xi \in \Irr(\FF G)$
such that, letting $i$ be a primitive
idempotent of $\FF K$ satisfying
$E(\eta) \cong \FF K i$, then
$E(\xi) \cong \FF G i \cong
{}_G \Ind {}_K (E(\eta))$. By Mackey
decomposition, the restriction
${}_K \Res {}_G (E(\xi))$ is a sum of
$G$-conjugates of $E(\eta)$. So the
condition $E(\xi) \cong {}_G \Ind {}_K
(E(\eta))$ characterizes a bijection
$[\eta]_G \leftrightarrow \xi$.

Suppose $[\eta]_G \leftrightarrow \xi$.
Since $V(\eta)$ and $V(\xi)$, respectively,
are the unique simple modules of $\FF K$
and $\FF G$ not annihilated by $i$, Clifford's
Theorem implies that ${}_K \Res {}_G
(V(\xi))$ is a direct sum of $G$-conjugates
of $V(\eta)$. Since $i$ is primitive in
$\FF G$, we have $\dim_\FF(i V(\xi)) = 1$,
so each $G$-conjugate of $V(\xi)$ occurs
in ${}_K \Res {}_G (V(\eta))$ with
multiplicity $1$.
\end{proof}

We shall be applying the lemma in the
following special case. Let $S$ be a
finite $p$-group acting as automorphisms
on a finite group $K$. Via the canonical
isomorphism $K \cong K {\rtimes} 1$, we
embed $K$ in the semidirect product
$K {\rtimes} S$. Given $\eta \in \Irr(\FF K)$,
we write $\eta_{\rtimes S}$ to denote the
irreducible $\FF (K {\rtimes} S)$-character
such that $\eta_{\rtimes S}$ corresponds
to the $S$-orbit of $\eta$.

Let $P {\ua} Q_\phi$ be a generalized
piece of $\FF G$. Now $P {\ua} Q_\phi$
determines $Q_\phi$ only up to $P$-conjugacy,
but let us make a choice of $Q_\phi$ and write
$\oo{P} = N_P(Q)/Q$. Via the conjugation
action of $N_P(Q)$ on $C_G(Q)$, we allow
$\oo{P}$ to act as automorphisms on
$C_G(Q)$ and we form the semidirect
product $C_G(Q) {\rtimes} \oo{P}$. Via the
canonical isomorphism $\oo{P} \cong
1 {\rtimes} \oo{P}$, we embed $\oo{P}$ in
$C_G(Q) {\rtimes} \oo{P}$. We call
$P {\ua} Q_\phi$ {\bf substantive} when the
simple $\FF(C_G(Q) {\rtimes} \oo{P})$-module
$V(\phi_{\rtimes \oo{P}})$ restricts to an
$\FF \oo{P}$-module with a nonzero free
direct summand. That condition is clearly
independent of the choice of $Q_\phi$.
Observe that any piece of $\FF G$ is a
substantive generalized piece of $\FF G$.

The following criterion for substantivity
will be proved, in a stronger quantitative
form, in Proposition \ref{3.2} below.

\begin{pro} \label{2.3}
Let $P {\ua} Q_\phi$ be a generalized piece
of $\FF G$. Write $\ooN_P(Q) = N_P(Q)/Q$.
Write $N_P(Q_\phi)$ for the
stabilizer of $Q_\phi$ in $P$ and write
$\ooN_P(Q_\phi) = N_P(Q_\phi)/Q$. Let
$\ooN_P(Q_\phi) \leq S \leq \ooN_P(Q)$. Then
$P {\ua} Q_\phi$ is substantive if and only if
the simple $\FF(C_G(Q) {\rtimes} S)$-module
$V(\phi_{\rtimes \oo{S}})$ restricts to an
$\FF S$-module with a nonzero free direct
summand.
\end{pro}

The proposition gives the following
means of determining whether a given
generalized piece $P {\ua} Q_\phi$ is
substantive. Let $T = \ooN_P(Q_\phi)$.
In the evident way, $\FF C_G(Q)$ becomes
a $T$-algebra. Via the representation of
$V(\phi)$, we regard $\End_\FF(V(\phi))$
as a $T$-algebra. Since the cohomology
groups $H^1(S, \FF^\times)$ and
$H^2(S, \FF^\times)$ are trivial, the
$T$-algebra structure of $\End_\FF(V(\phi))$
enriches, in a unique way, to an interior
$T$-algebra structure. Thus, $V(\phi)$
becomes an $\FF T$-module. The
generalized piece $P {\ua} Q_\phi$ is
substantive if and only if, regarding
$V(\phi)$ as an $\FF T$-module, the
regular $\FF T$-module occurs as a
direct summand of $V(\phi)$.

We note an immediate corollary of
the proposition.

\begin{cor} \label{2.4}
Let $P {\ua} Q_\phi$ be a generalized
piece of $\FF G$ and let $Q \leq P' \leq P$.

\tinyskip
\noin {\bf (1)} If $P {\ua} Q_\phi$ is
substantive, then $P' {\ua} Q_\phi$ is
substantive.

\tinyskip
\noin {\bf (2)} If $N_P(Q_\phi) \leq P'$
and $P' {\ua} Q_\phi$ is substantive,
then $P {\ua} Q_\phi$ is substantive.
\end{cor}

For any piece $P {\ua} Q_\phi$ of $\FF G$, we
now construct an $\FF(G {\times} P)$-module
$\Dia_G(P {\ua} Q_\phi)$, called the
{\bf diagonal module} of $P {\ua} Q_\phi$.
Still working with a choice of $Q_\phi$, let
$$N = N_{G \times P}(\Delta(Q)) =
  (C_G(Q) {\times} 1) \Delta(N_P(Q))
  \; , \;\;\;\; \;\;\;\;
  \ooN = \ooN_{G \times P}(\Delta(Q))
  \cong C_G(Q) {\rtimes} \oo{P} \; .$$
Via the canonical isomorphism
$C_G(Q) \cong (C_G(Q) {\times} 1)
\Delta(Q)/\Delta(Q)$, we embed
$C_G(Q)$ in $\ooN$ and form the
induced $\FF \ooN$-module
$$\oo{\Dia}_G^0(P {\ua} Q_\phi) =
  {}_{\ooN} \Ind {}_{C_G(Q)} (E(\phi))
  \cong {}_{\ooN} \Iso
  {}_{C_G(Q) \rtimes \oo{P}}
  (E(\phi_{\rtimes \oo{P}}))$$
which is indecompoable and projective.
We define the inflated $\FF N$-module
$$\Dia_G^0(P {\ua} Q_\phi) = {}_N \Inf
  {}_{\ooN} (\oo{\Dia}_G^0(P {\ua} Q_\phi))$$
which is indecomposable with vertex
$\Delta(Q)$. We define
$\Dia_G(P {\ua} Q_\phi)$ to be the
indecomposable $\FF(G {\times} P)$-module
with vertex $\Delta(Q)$ in Green
correspondence with $\Dia_G^0
(P {\ua} Q_\phi)$. It is easy to check that,
given a $P$-conjugate $Q'_{\phi'}$ of
$Q_\phi$, then $\Dia_G(P {\ua} Q_\phi)
\cong \Dia_G(P {\ua} Q'_{\phi'})$. Thus,
$\Dia_G(P {\ua} Q_\phi)$ is determined
by $P {\ua} Q_\phi$ up to isomorphism,
independently of the choice of $Q_\phi$.
The next result tells us that, in fact,
$\Dia_G(P {\ua} Q_\phi)$ is uniquely
determined by $P {\ua} Q_\phi$.

\begin{lem} \label{2.5}
Let $P$ be a $p$-subgroup of $G$ and let
$Q_\phi$ and $Q'_{\phi'}$ be pieces of
$\FF G$ such that $P \geq Q$ and
$P \geq Q'$. Then $\Dia_G(P {\ua} Q_\phi)
\cong \Dia_G(P {\ua} Q'_{\phi'})$ if and only
if $P {\ua} Q_\phi = P {\ua} Q'_{\phi'}$.
\end{lem}

\begin{proof}
The conclusion in one direction has been
observed already. For the converse,
suppose $\Dia_G(P {\ua} Q_\phi) \cong
\Dia_G(P {\ua} Q'_{\phi'})$. By considering
vertices, $\Delta(Q)$ and $\Delta(Q')$ are
$G {\times} P$-conjugate, hence $Q$ and
$Q'$ are $P$-conjugate and we may assume
that $Q = Q'$. By considering Green
correspondents, $E(\phi_{\rtimes \oo{P}})
\cong E(\phi'_{\rtimes \oo{P}})$. By Lemma
\ref{2.2}, $\phi$ and $\phi'$ lie in the same
$\oo{P}$-orbit of $\Irr(\FF C_G(Q))$, hence
$Q_\phi$ and $Q_{\phi'}$ are $P$-conjugate.
\end{proof}

We shall be making use of the following
characterization of the diagonal
module $\Dia(P {\ua} Q_\phi)$.

\begin{lem} \label{2.6}
Given a generalized piece $P {\ua} Q_\phi$
on $\FF G$, then
$$\Dia_G(P {\ua} Q_\phi) \cong {}_{G \times P}
  \Ind {}_{G \times Q} (\Dia_G(Q_\phi)) \; .$$
\end{lem}

\begin{proof}
By its definition, $\Dia_G(P {\ua} Q_\phi)$ is
the isomorphically unique indecomposable
$\FF(G {\times} P)$-module with vertex
$\Delta(Q)$ that appears as a direct
summand of the $\FF(G {\times} P)$-module
$$L = {}_{G \times P} \Ind
  {}_{N_{G \times P}(\Delta(Q))} \Inf
  {}_{\ooN_{G \times P}(\Delta(Q))} \Ind
  {}_{C_G(Q)} (E(\phi)) \; .$$
As a special case, $\Dia_G(Q_\phi)$ is the
isomorphically unique indecomposable
$\FF(G {\times} Q)$-module with vertex
$\Delta(Q)$ that appears as a direct
summand of
$$M = {}_{G \times Q} \Ind
  {}_{N_{G \times Q}(\Delta(Q))} \Inf
  {}_{C_G(Q)} (E(Q))$$
where the inflation is via the canonical
epimorphism $N_{G \times Q}(\Delta(Q))
\rightarrow \ooN_{G \times Q}(\Delta(Q))
\cong C_G(Q)$. Using the Mackey formula
for bisets in Bouc \cite[2.3.24]{Bou10}, we
obtain an equality of bisets
$${}_{N_{G \times P}(\Delta(Q))} \Inf
  {}_{\ooN_{G \times P}(\Delta(Q))}
  \Ind {}_{C_G(Q)} \cong
  {}_{N_{G \times P}(\Delta(Q))} \Ind
  {}_{N_{G \times Q}(\Delta(Q))}
  \Inf {}_{C_G(Q)} \; .$$
Hence, $L \cong {}_{G \times P} \Ind
{}_{G \times Q} (M)$. Therefore, writing
$D = {}_{G \times P} \Ind {}_{G \times Q}
(\Dia_G(Q_\phi))$, then $D$ is a direct
summand of $L$. By Green's
Indecomposability Criterion, $D$ is
indecomposable with vertex $\Delta(Q)$. By
the uniqueness of $\Dia_G(P {\ua} Q_\phi)$,
we have $\Dia_G(P {\ua} Q_\phi) \cong D$.
\end{proof}

When we write an $\FF(G {\times} P)$-module
in the form ${}_G M {}_P$, we are indicating
that the actions of $G {\times} 1$ and
$1 {\times} P$ are by left translation and
right translation, respectively. Given
$\FF G$-modules $L$ and $M$, we write
$L \smid M$ when $L$ is isomorphic to
a direct summand of $M$. A stronger
quantitative version of the next result will
appear below as Proposition \ref{3.5}

\begin{pro} \label{2.7}
A generalized piece $P {\ua} Q_\phi$ on
$\FF G$ is substantive if and only if
$$\Dia_G(P {\ua} Q_\phi) \smid
  {}_G \FF G {}_P \; .$$
\end{pro}

Our main result is the following classification
of the pointed $p$-groups on $\FF G$. We
shall prove it in the next section.

\begin{thm} \label{2.8}
There is a $G$-equivariant bijective
correspondence between:

\tinyskip
\noin $\bullet$ the pointed $p$-groups
$P_\alpha$ on $\FF G$,

\tinyskip
\noin $\bullet$ the substantive
generalized pieces $P' {\ua} Q_\phi$,

\tinyskip
\noin such that $P_\alpha \leftrightarrow
P' {\ua} Q_\phi$ if and only if $P = P'$ and,
letting $Q_\delta$ be the local pointed group
on $\FF G$ corresponding to the piece
$Q_\phi$, also letting $i \in \alpha$, the
following two equivalent conditions hold:

\tinyskip
\noin {\bf (a)} $Q_\delta$ is a maximal
local pointed subgroup of $P_\alpha$,

\tinyskip
\noin {\bf (b)} we have $\Dia_G(P {\ua}
Q_\phi) \cong {}_G (\FF G i) {}_P$.
\end{thm}

The bijective correspondences in
Theorems \ref{2.1} and \ref{2.8} are
compatible with blocks in the following
ways. Let $b$ be a block of $\FF G$. A
piece $P_\theta$ of $\FF G$ is called a
{\bf piece} of $\FF G b$ provided
$\br_P(b)$ acts as the identity on the
$\FF C_G(P)$-module $E(\theta)$.
Obviously, letting $P_\gamma$ be the
local pointed group on $\FF G$
corresponding to $P_\theta$, then
$P_\theta$ is a piece of $\FF G b$
if and only if $P_\gamma$ is a local
pointed group on $\FF G b$. A
generalized piece $P {\ua} Q_\phi$
of $\FF G$ is called a {\bf generalized
piece} of $\FF G b$ provided $Q_\phi$
is a piece of $\FF G b$. Theorem \ref{2.8}
has the following immediate corollary.

\begin{cor} \label{2.9}
Let $b$ be a block of $\FF G$. Let
$P_\alpha$ be a pointed $p$-group on
$\FF G$. Let $P {\ua} Q_\phi$ be the
substantive generalized piece of $\FF G$
corresponding to $P_\alpha$. Then the
following three conditions are equivalent:

\tinyskip
\noin {\bf (a)} $P_\alpha$ is a pointed
$p$-group on $\FF G b$,

\tinyskip
\noin {\bf (b)} $P {\ua} Q_\phi$ is
a piece of $\FF G b$,

\tinyskip
\noin {\bf (c)} we have
$\Dia_G(P {\ua} Q_\phi) \smid
{}_G (\FF G b) {}_P$.
\end{cor}

\section{Stronger quantitative results}
\label{3}

We shall define the absolute multiplicity
of a generalized piece. That will enable
us to prove stronger quantitative versions
of Propositions \ref{2.3}, \ref{2.7} and
Theorem \ref{2.8}.

For any pointed group $U_\mu$, we write
$m(U_\mu)$ to denote the absolute
multiplicity of $U_\mu$, we mean to say,
the maximal size of a set of mutually
orthogonal elements of $\mu$. Given a
piece $P_\theta$ of $\FF G$, we define
the {\bf absolute multiplicity} of $P_\theta$
to be
$$m(P_\theta) = \theta(1) =
  \dim_\FF(V(\theta)) \; .$$
The next remark says that the bijective
correspondence in Theorem 2.1 preserves
absolute multiplicities.

\begin{rem} \label{3.1}
Given a local pointed group $P_\gamma$
on $\FF G$ with corresponding piece
$P_\theta$ on $\FF G$, then
$m(P_\gamma) = m(P_\theta)$.
\end{rem}

\begin{proof}
Letting $i \in \gamma$, then $m(P_\gamma)$
and $m(P_\theta)$ are both equal to the
multiplicity of the projective indecomposable
$\FF C_G(P)$-module $\FF C_G(P) \br_P(i)$
as a direct summand of the regular
$\FF C_G(P)$-module.
\end{proof}

To prove the results in the previous section,
we shall need to extend the notion of
absolute multiplicity to generalized
pieces. Given $\FF G$-modules $L$ and
$M$ with $L$ indecomposable, we write
$m(L, M)$ to denote the multiplicity of
$L$ as a direct summand of $M$.

Let $P {\ua} Q_\phi$ be a
generalized piece on $\FF G$. As before,
we make a choice of $Q_\phi$. Again, we
write $\oo{P} = N_P(Q)$ and we consider the
simple $\FF(C_G(Q) {\rtimes} \oo{P})$-module
$V(\phi_{\rtimes \oo{P}})$. We define the
{\bf absolute multiplicity} of $P {\ua} Q_\phi$
to be the natural number
$$m(P {\ua} Q_\phi) = m(\FF \oo{P},
  {}_{\oo{P}} \Res {}_{C_G(Q) \rtimes \oo{P}}
  (V(\phi_{\rtimes \oo{P}})))$$
where $\FF \oo{P}$ denotes the regular
$\FF \oo{P}$-module.
Plainly, $m(P {\ua} Q_\phi)$ is well-defined,
independently of the choice of $Q_\phi$.
Observe that, for a piece $P_\theta$ of
$\FF G$, the absolute multiplicity
$m(P_\theta) = m(P {\ua} P_\theta)$ is
unambiguous. The generalized piece
$P {\ua} Q_\phi$ is substantive if and only
if $m(P_\theta) \neq 0$.

The next result is a stronger quantitative
version of Proposition \ref{2.3}.

\begin{pro} \label{3.2}
Let $P {\ua} Q_\phi$ be a generalized piece
on $\FF G$. Let $\ooN_P(Q_\phi) \leq S \leq
\ooN_P(Q)$. Then
$$m(P {\ua} Q_\phi) =
  m(\FF S, {}_S \Res {}_{C_G(Q) \rtimes S}
  (V(\phi_{\rtimes S}))) \; .$$
\end{pro}

\begin{proof}
Write $C = C_G(Q)$ and $T =
\ooN_P(Q_\phi)$. Define $m_S =
m(\FF S, {}_S \Res {}_{C_G(Q) \rtimes S}
(V(\phi_{\rtimes S})))$. By considering the
case where $S = \ooN_P(Q)$, we see that
it suffices to show that $m_T = m_S$. We
have ${}_C \Res {}_{C \rtimes T}
(V(\phi_{\rtimes T})) \cong V(\phi)$ and
$${}_C \Res {}_{C \rtimes S}
  (V(\phi_{\rtimes S})) \cong
  \bigoplus_{sT \subseteq S}
  {}^s V(\phi) \; .$$
So $\dim_\FF(V(\phi_{\rtimes S})) =
|S : T| \, \dim_\FF(V(\phi_{\rtimes T}))$.

Since $C {\rtimes} T$ is subnormal in
$C {\rtimes} S$, Clifford's Theorem implies
that ${}_{C \rtimes T} \Res {}_{C \rtimes S}
(V(\phi_{\rtimes S}))$ is semisimple. But
$V(\phi)$ occurs in the semisimple
$\FF C$-module ${}_C \Res
{}_{C \rtimes S} (V(\phi_{\rtimes S}))$
and $V(\phi_{\rtimes T})$ is the
isomorphically unique simple
$\FF(C {\rtimes} T)$-module such
that $V(\phi)$ occurs in the semisimple
$\FF(C {\rtimes} T)$-module ${}_C \Res
{}_{C \rtimes T} (V(\phi_{\rtimes T}))$.
Therefore, $V(\phi_{\rtimes T})$ occurs in
the semisimple $\FF(C {\rtimes} T)$-module
${}_{C \rtimes T} \Res {}_{C \rtimes S}
(V(\phi_{\rtimes S}))$. By Frobenius
reciprocity, $V(\phi_{\rtimes S})$ is
isomorphic to a submodule of
${}_{C \rtimes S} \Ind {}_{C \rtimes T}
(V(\phi_{\rtimes T}))$. A consideration
of dimensions yields
$$V(\phi_{\rtimes S}) \cong
  {}_{C \rtimes S} \Ind {}_{C \rtimes T}
  (V(\phi_{\rtimes T})) \; .$$
By Mackey decomposition,
$${}_S \Res {}_{C \rtimes S}
  (V(\phi_{\rtimes S})) \cong
  {}_S \Ind {}_T \Res {}_{C \rtimes T}
  (V(\phi_{\rtimes T})) \; .$$
The required equality $m_T = m_S$ follows
because, given any indecomposable direct
summand $M$ of ${}_T \Res {}_{C \rtimes T}
(V(\phi_{\rtimes T}))$, then $M$ is free
if and only if ${}_S \Ind {}_T (M)$ is free.
\end{proof}

We shall be needing the following result
of Brou\'{e} \cite[3.2]{Bro85}. Given a
$p$-subgroup $S$ of $G$ and an
$\FF G$-module $M$, we define the
$\FF \ooN_G(S)$-module $M(S)$ to be the
quotient of $M^S$ by the sum of the images of the
trace maps $\tr_T^S : M^T \rightarrow M^S$,
running over the strict subgroups $T < S$.

\begin{pro} \label{3.3}
{\rm (Brou\'{e}.)} Let $S$ be a $p$-subgroup
of $G$, let $E$ be an indecomposable projective
$\FF \ooN_G(S)$-module, and let $F$ be the
indecomposable $\FF G$-module with vertex
$S$ in Green correspondence with the inflated
$\FF N_G(S)$-module ${}_{N_G(S)} \Inf
{}_{\ooN_G(S)} (E)$. Let $M$ be a
$p$-permutation $\FF G$-module. Then
$m(F, M) = m(E, M(S))$.
\end{pro}

Another necessary ingredient is the following
result of Robinson \cite[Proposition 1]{Rob89}.

\begin{pro} \label{3.4}
{\rm (Robinson.)} Given $\xi \in \Irr(\FF G)$
and a $p$-subgroup $S$ of $G$, then
$$m(E(\xi), \FF G/S) =
m(\FF S, {}_S \Res {}_G (V(\xi))) \; .$$
\end{pro}

The next result implies Proposition \ref{2.7}
and, more precisely, it characterises the
multiplicity of a generalized piece in terms
of the associated diagonal module. 

\begin{pro} \label{3.5}
Given a generalized piece $P {\ua} Q_\phi$
on $\FF G$, then
$$m(P {\ua} Q_\phi) =
  m(\Dia_G(P {\ua} Q), {}_G \FF G {}_P) \; .$$
\end{pro}

\begin{proof} Let $m = m(\Dia_G(P {\ua} Q),
{}_G \FF G {}_P)$. By Proposition \ref{3.3},
$$m(\Dia_G(P {\ua} Q_\phi), M) =
  m(\oo{\Dia}_G^0(P {\ua} Q_\phi),
  M(\Delta(Q)))$$
for any $p$-permutation
$\FF(G {\times} P)$-module $M$. Put
$M = {}_G \FF G {}_P$. Let $N$ and $\ooN$
be as in Section \ref{2}. Define
$L = {}_N \Res {}_{G \times P} (M)$. Then
$M(\Delta(Q)) = L(\Delta(Q))$ and
$$m = m(\oo{\Dia}_G^0(P {\ua} Q_\phi),
  L(\Delta(Q))) \; .$$
Let $\Gamma \subseteq G$ such that
$1 \in \Gamma$ and $\{ (g, 1) :
g \in \Gamma \}$ is a set of representatives
of the double cosets $N \back (G {\times} P) /
\Delta(P)$. For each $g \in \Gamma$, we
define a permutation $\FF N$-module
$L_g = \FF N / H_g$ where $H_g =
N \cap {}^g \Delta(P)$. Since $M \cong
\FF(G \times P) / \Delta(P)$, Mackey
decomposition yields
$$L(\Delta(Q)) \cong \bigoplus_{g \in \Gamma}
  L_g(\Delta(Q)) \; .$$

Fix $g \in \Gamma$ such that $L_g(\Delta(Q))
\neq 0$. Then $\Delta(Q) \leq H_g \leq
{}^g \Delta(P)$, hence $g \in C_G(Q)$ and
$(g, 1) \in N$. The condition $1 \in \Gamma$
implies that $g = 1$. Therefore, $L(\Delta(Q))
\cong L_1(\Delta(Q))$. Since $H_1 =
\Delta(N_P(Q))$, we have $L(\Delta(Q))
\cong \FF \ooN/\oo{P}$. We have shown that
$$m = m({}_{\ooN} \Ind {}_{C_G(Q)}
  (E(\phi)), \FF \ooN / \oo{P}) \; .$$
In view of the isomorphism $\ooN \cong
C_G(Q) {\rtimes} \oo{P}$, Lemma \ref{2.2}
yields
$$m = m(E(\phi_{\rtimes \oo{P}}), \FF
  (C_G(Q) {\rtimes} \oo{P}) / \oo{P}) \; .$$
Lemma \ref{3.4} now implies that
$m = m(P {\ua} Q_\phi)$.
\end{proof}

We shall be needing two abstract lemmas.
We write $\leq$ to denote the usual partial
ordering on the idempotents of a ring.

\begin{lem} \label{3.6}
For any local pointed group $P_\gamma$
on $\FF G$, there exists $i \in \gamma$
such that $i \leq \br_P(i)$.
\end{lem}

\begin{proof}
Let $i_0 \in \gamma$. Then
$\br_P(i_0 \br_P(i_0)) = \br_P(i_0)
\not\in J(\FF C_G(P))$. So
$i_0 \br_P(i_0) \not\in J((\FF G)^P)$. It
follows that $i \leq \br_P(i_0)$ for some
$i \in \gamma$. We have $\br_P(i) \leq
\br_P(\br_P(i_0)) = \br_P(i_0)$. But
$\br_P(i_0)$ is a primitive idempotent
of $\FF C_G(P)$. Therefore,
$\br_P(i) = \br_P(i_0)$.
\end{proof}

We point out that the proof of the
lemma yields a stronger result, namely,
that for all $i_0 \in \gamma$, there
exists $i \in \gamma$ satisfying
$i \leq \br_P(i) = \br_P(i_0)$.

\begin{lem} \label{3.7}
Let $P_\alpha$ be a pointed $p$-group on
$\FF G$. Given $i \in \alpha$, then the
$\FF(G {\times} P)$-module
${}_G (\FF G i) {}_P$ is indecomposable.
Given $Q \leq P$, then $Q$ is a defect
group of $P_\alpha$ if and only if $\Delta(Q)$
is a vertex of ${}_G (\FF G i) {}_P$.
\end{lem}

\begin{proof}
Writing $\circ$ to indicate an
opposite algebra, there is an interior
$P$-algebra isomorphism
$$\End_{\FF(G \times 1)} (\FF G i) \cong
  (i \FF G i)^\circ$$
such that, given $r \in \End_{\FF(G \times 1)}
(\FF G i)$ and $a \in i \FF G i$, then
$r \leftrightarrow a^\circ$ provided
$r(x) = xa$ for all $x \in \FF G i$. Hence
$\End_{\FF(G \times P)} (\FF G i) \cong
((i \FF G i)^P)^\circ$, which is a local algebra.
Therefore, ${}_G (\FF G i) {}_P$ is
indecomposable. Since ${}_G (\FF G i) {}_P$
is a direct summand of the
$\FF(G {\times} P)$-module ${}_G \FF G {}_P
\cong \FF(G {\times} P)/\Delta(P)$, some
vertex of ${}_G (\FF G i) {}_P$ is contained
in $\Delta(P)$.

Suppose $Q$ is a defect group of
$P_\alpha$. Let $a \in (i \FF G i)^Q$ such
that $i = \tr_Q^P(a)$. Let $r \in 
\End_{\FF(G \times Q)}(\FF G i)$ such
that $r \leftrightarrow a^\circ$. Then
$\tr_{G \times Q}^{G \times P}(r) =
\tr_{\Delta(Q)}^{\Delta(P)}(r) =
\id_{\FF G i}$. So $G {\times} Q$ contains a
vertex $S$ of $\FF G i$. But
$S \leq {}^{(g, u)} \Delta(P)$
for some $(g, u) \in G {\times} P$. So the
vertex ${}^{(u g^{-1}, 1)} S$ of $\FF G i$
is contained in the subgroup $(G {\times} Q)
\cap \Delta(P) = \Delta(Q)$.

For the reverse inclusion, suppose
$\Delta(Q)$ is a vertex of $\FF G i$. Let
$r \in \End_{\FF(G \times Q)}(\FF G i)$ such
that $\tr_{G \times Q}^{G \times P}(r)
= \id_{\FF G i}$. Let $a \in (i \FF G i)^Q$
such that $r \leftrightarrow a^\circ$. Then
$i = \tr_Q^P(a)$. We deduce that $Q$
contains a defect group of $P_\alpha$.
\end{proof}

\begin{pro} \label{3.8}
Let $P_\gamma$ be a local pointed group
on $\FF G$. Let $P_\theta$ be the piece of
$\FF G$ corresponding to $P_\gamma$.
Let $i \in \gamma$. Then $\Dia_G(P_\theta)
\cong {}_G (\FF G i) {}_P$.
\end{pro}

\begin{proof}
It is easy to check that the isomorphism class
of the $\FF(G {\times} P)$-module
${}_G (\FF G i) {}_P$ is determined by
$P_\gamma$, independently of $i$. So, in
view of Lemma \ref{3.6}, we may assume
that $i \leq \br_P(i)$. Let
$\br_P(i) = i + \sum_k k$
be a primitive idempotent decomposition
in $(\FF G)^P$. By Lemma \ref{3.7},
$${}_G (\FF G \br_P(i)) {}_P =
  {}_G (\FF G i) {}_P \oplus \bigoplus_k
  {}_G (\FF G k) {}_P$$
as a direct sum of indecomposable
$\FF(G {\times} P)$-modules. Each $k$
has a defect group $Q_k$ strictly
contained in $P$. By the same lemma,
${}_G (\FF G i) {}_P$ has vertex $\Delta(P)$
while each ${}_G (\FF G k) {}_P$ has vertex
$\Delta(Q_k)$. Inflating via the canonical
epimorphism $N_{G \times P}(\Delta(P))
\rightarrow C_G(P)$, we have
$$\FF C_G(P) \br_P(i) \cong
  {}_{N_{G \times P}(\Delta(P)} \Inf
  {}_{C_G(P)} (E(\theta)) \cong
  \Dia_G^0 (P_\theta) \; .$$
As a direct sum of $\FF N_{G \times P}
(\Delta(P))$-modules,
$$\FF G \br_P(i) = \FF C_G(P) \br_P(i)
  \oplus \FF(G - C_G(P)) \br_P(i) \; .$$
The conjugation action of $P$ on $G - C_G(P)$
has no fixed points, so each indecomposable
direct summand of $\FF(G - C_G(P)) \br_P(i)$
has a vertex strictly contained in $\Delta(P)$.
The required isomorphism now follows from
the Green Correspondence Theorem.
\end{proof}

The next result is a theorem of Puig that can
be found in Linckelmann \cite[5.12.20]{Lin18}
or Th\'{e}venaz \cite[18.3]{The95}.

\begin{thm} \label{3.9}
{\rm (Puig.)} Let $Q_\delta$ be a local
pointed group on a $G$-algebra $A$ over
$\FF$. Let $P_\alpha$ be be a pointed
group on $A$ such that $Q \leq P$. Let
$i \in \alpha$. Then $Q_\delta$ is a
maximal local pointed subgroup of
$P_\alpha$ if and only if there exists
$j \in \delta$ such that $i = \tr_Q^P(j)$
and $\{ {}^x j : xQ \subseteq P \}$ is a
set of mutually orthogonal idempotents.
\end{thm}

We shall also be using the following
result of \cite[3.1]{BG22}.

\begin{thm} \label{3.10}
Let $Q_\delta$ be a local pointed group
on a $G$-algebra $A$ over $\FF$. Let
$P$ be a $p$-subgroup of $G$ such that
$Q \leq P$. Then there exists at most one
point $\alpha$ of $P$ on $A$ such that
$Q_\delta$ is a maximal local pointed
subgroup of $P_\alpha$.
\end{thm}

\begin{lem} \label{3.11}
Let $P$ be a $p$-subgroup of $G$,
let $\alpha$ and $\alpha'$ be points
of $P$ on $\FF G$ and let $i \in \alpha$
and $i' \in \alpha'$. Then ${}_G (\FF G i)
{}_P \cong {}_G (\FF G i') {}_P$ if and
only if $\alpha = \alpha'$.
\end{lem}

\begin{proof}
If $\alpha = \alpha'$, then the isomorphism
is clear. Conversely, assume the
isomorphism. Lemma \ref{3.7} implies
that $P_\alpha$ and $P_{\alpha'}$ have
a common defect group $Q$. Let $\delta$
and $\delta'$ be local points of $Q$ on
$\FF G$ such that $Q_\delta$ and
$Q_{\delta'}$ are maximal local pointed
subgroups of $P_\alpha$ and
$P_{\alpha'}$, respectively. By
Theorem \ref{3.9}, there exist
$j \in \delta$ and $j' \in \delta'$ such
that $i = \tr_Q(j)$ as a sum of
$|P : Q|$ mutually orthogonal
idempotents and similarly for $i'$
and $j'$. We have
$${}_G (\FF G i) {}_P \cong
  {}_{G \times P} \Ind {}_{G \times Q}
  (\FF G j) \; .$$
Restricting to $G {\times} Q$ and then
applying Mackey decomposition,
$$(\FF G i)(\Delta(Q)) \cong
  ({}_{G \times Q} \Res {}_{G \times P}
  \Ind {}_{G \times Q} (\FF G j))
  (\Delta(Q)) \cong
  \bigoplus_{xQ \subseteq N_P(Q)}
  \FF C_G(Q) \, \br_Q({}^x j)$$
and similarly for $i'$ and $j'$. Therefore,
$\FF C_G(Q) \, \br_Q(j') \cong
\FF C_G(Q) \, \br_Q({}^x j)$ for some
$x \in N_P(Q)$. Hence, $\delta' =
{}^x \delta$. So we may assume that
$\delta' = \delta$. By Theorem
\ref{3.10}, $\alpha = \alpha'$.
\end{proof}

We now prove Theorem \ref{2.8}. For each
$p$-subgroup $P$ of $G$, Lemmas
\ref{2.5} and \ref{3.11} imply that the
isomorphism classes of $\Dia_G(P {\ua} Q_\phi)$
and ${}_G (\FF G i) {}_P$ determine
$P {\ua} Q_\phi$ and $P_\alpha$, respectively.
So Propositions \ref{3.5} and \ref{3.8} imply
that condition (b) characterizes a bijective
correspondence $P_\alpha \leftrightarrow
P {\ua} Q_\phi$.

It remains only to show that
the conditions (a) and (b) are equivalent. Let
$Q_\delta$ be the local pointed group
corresponding to $Q_\phi$. Let $j \in \delta$.
By Proposition \ref{3.8}, $\Dia_G(Q_\phi) \cong
{}_G (\FF G j) {}_Q$. So by Lemma \ref{2.6},
$$\Dia_G(P {\ua} Q_\phi) \cong {}_{G \times P}
  \Ind {}_{G \times Q} (\FF G j) \; .$$
If (a) holds then, by Theorem \ref{3.9}, we
can take the choice of $j$ to be such that
$i = \tr_Q^P(j)$ as a sum of $|P : Q|$ mutually
orthogonal idempotents. We deduce (b).

Conversely, assume (b). By Lemma \ref{3.7},
$P_\gamma$ has defect group $Q$. Let
$\delta'$ be a point of $Q$ on $\FF G$ such
that $P_\gamma$ has maximal local pointed
subgroup $Q_{\delta'}$. By what we have
already shown, $\Dia_G(P {\ua} Q_{\delta'})
\cong \Dia_G(P {\ua} Q_\delta)$. By Lemma
\ref{2.5}, $P {\ua} Q_{\delta'} =
P {\ua} Q_\delta$. So $Q_{\delta'}$ and
$Q_\delta$ are $P$-conjugate, and (a) holds.
The proof of Theorem \ref{2.8} is complete.

Let us note that the bijective correspondence
in that theorem preserves absolute multiplicities
in the following sense.

\begin{pro} \label{3.12}
Let $P_\alpha$ be a pointed $p$-group on
$\FF G$. Let $P {\ua} Q_\phi$ be the substantive
generalized piece on $\FF G$ corresponding to
$P_\alpha$. Then $m(P_\alpha) =
m(P {\ua} Q_\phi)$.
\end{pro}

\begin{proof}
This follows immediately from
Proposition \ref{3.5}
\end{proof}

\section{The principal $2$-blocks of $S_4$
and $S_5$}

\label{4}
By a method that involves calculation
of the absolute multiplicities, we shall
determine the relative multiplicities between
the substantive generalized pieces for the
principal $2$-blocks of $S_4$ and $S_5$.
The method is based on the following
immediate implication of Theorem \ref{2.8},
Corollary \ref{2.9}, Proposition \ref{3.12}.

\begin{cor} \label{4.1}
Given a $p$-subgroup $P$ of $G$, then
we have an $\FF(G {\times} P)$-isomorphism
$${}_G \FF G {}_P \cong
  \bigoplus_{Q_\phi} m(P {\ua} Q_\phi)
  \, {}_G (\FF G i_{Q_\phi}) {}_P$$
where $Q_\phi$ runs over the
$P$-conjugacy classes of pieces of
$\FF G$ such that the generalized piece
$P {\ua} Q_\phi$ is substantive, and
$i_{Q_\phi}$ is an element of the point
$\alpha$ of $P$ on $\FF G$ such that
$P_\alpha$ is the local pointed group on
$\FF G$ corresponding to
$P {\ua} Q_\phi$. Furthermore, given a
block $b$ of $\FF G$, then
$${}_G (\FF G b) {}_P \cong
  \bigoplus_{Q_\phi} m(P {\ua} Q_\phi)
  \, {}_G (\FF G i_{Q_\phi}) {}_P$$
where $Q_\phi$ now runs over the
$P$-conjugacy classes of pieces of
$\FF G b$ such that $P {\ua} Q_\phi$ is
substantive.
\end{cor}

Consider a block $b$ of $\FF G$. Let us
make some comments on some combinatorial
structures possessed by the set $\cP(\FF G b)$
of pieces of $\FF G b$ and the set
$\cP^|(\FF G b)$ of substantive generalized
pieces of $\FF G b$. We understand a
{\bf multiposet} to be a poset such that
each inclusion $x' \leq x$ is associated
with a natural number $m(x', x)$ called
the {\bf multiplicity} of $x'$ in $x$. We
regard $\cP^|(\FF G b)$ as a poset
such that, given $P' {\ua} Q'_{\phi'},
P {\ua} Q_\phi \in \cP^|(\FF G b)$, then
$P' {\ua} Q'_{\phi'} \leq P {\ua} Q_\phi$
provided $P' \leq P$ and the relative
multiplicity $m(P' {\ua} Q'_{\phi'},
P {\ua} Q_\phi)$ is nonzero. We regard
$\cP^|(\FF G b)$ as a multiposet whose
multiplicities are the relative multiplicities.
In an evident sense, $\cP(\FF G b)$ is
a submultiposet of $\cP^|(\FF G b)$.

By the matrix relation for relative
multiplicities in Section \ref{1}, the whole
family of relative multiplicities
$m(P' {\ua} Q'_{\phi'}, P {\ua} Q_\phi)$ is
determined by those relative multiplicities
such that $P' \leq P$ and $|P : P'| = p$. The
Hasse diagram for $\cP^|(\FF G b)$ as a poset
has an upwards line from $P' {\ua} Q'_{\phi'}$
to $P {\ua} Q_\phi$ if and only if that
condition on $P'$ and $P$ holds. So the
structure of $\cP^|(\FF G b)$ is determined
by an enriched Hasse diagram where any
upwards line from an element
$P' {\ua} Q'_{\phi'}$ to an element
$P {\ua} Q_\phi$ is labled with
$m(P' {\ua} Q'_{\phi'}, P {\ua} Q_\phi)$.

Put $p = 2$. Let $H = S_4$. The principal
block of $\FF H$ is the unique block of
$\FF H$. Let $D$ be a Sylow $2$-subgroup
of $H$. We have $D \cong D_8$, the dihedral
group of order $8$. Let
$$\SS = \{1, C_2, C'_2, Z,
  C_4, V_4, V'_4, D \}$$
be a set of representatives of the
$D$-conjugacy classes of subgroups
of $D$, named according to the isomorphism
classes and such that $Z =_G C_2 < V_4$.
Let $\PP_H$ be the submultiposet of
$\cP^|(\FF H)$ consisting of those
substantive generalized pieces on $\FF H$
that have the form $P {\ua} Q_\phi$ where
$P, Q \in \SS$. The elements of $\PP_H$
comprise a set of representatives for the
$D$-conjugacy classes of substantive
generalized pieces of $\FF H$ having the
form $P {\ua} Q_\phi$ where $P \leq D$. So
a Hasse diagram for $\PP_H$, labelled with
the relative multiplicities, will supply a
complete description of the multiposet
$\cP^|(\FF H)$.

Determining the pieces in $\PP_H$ will
be straightforward. To find the other
elements of $\PP_H$, we shall use the
Clifford-theoretic criterion for substantivity
of a generalized piece. That criterion was
presented above immediately following the
statement of Proposition \ref{2.3}.

Let us set up some conventions of notation
for expressing pieces of $\FF H$ concisely.
For any subgroup $L \leq H$, we write
$\Irr(\FF L) = \{ \theta_1^L, \theta_2^L, ... \}$,
enumerated such that $\theta_1^L$ is the
trivial $\FF L$-character. For any $p$-subgroup
$P \leq H$, we write $P_i = P_{\theta_i^C}$
where $C = C_H(P)$. We have $\Irr(\FF H) =
\{ \theta_1^H, \theta_2^H \}$. For all
$P \in \SS - \{ 1 \}$, the centralizer $C_H(P)$
is a $p$-group, so $\Irr(\FF C_G(P)) =
\{ \theta_1^{C_H(P)} \}$. Therefore, the
pieces in $\PP_H$ are
$$1_1, 1_2, (C_2)_1, (C'_2)_1, Z_1,
  (C_4)_1, (V_4)_1, (V'_4)_1, D_1 .$$

We claim that the only other element of
$\PP_H$ is $C'_2 {\ua} 1_2$. To demonstrate
the claim, we first note that, for all $Q \in \SS$,
we have $\theta_1^{C_H(Q)}(1) = 1$, so there
is no element $P {\ua} Q_1 \in \PP_H$
with $P > Q$. The kernel of $V(\theta_2^H)$
is $V_4$. So, given $P \in \SS$, then
${}_P \Res {}_H (V(\theta_2^H))$ has a
nonzero free direct summand if and only if
$P \in \{ \{ 1, C'_2 \}$. The claim is
established.

Recall, the absolute multiplicity of a piece
is the degree of the associated irreducible
character. For the sole element of
$\PP_H$ that is not a piece, we have
${}_{C'_2} \Res {}_H (V(\theta_2^H))
\cong \FF C'_2$ and the absolute
multiplicity is $m(C'_2 {\ua} 1_2) = 1$.
Applying Corollary \ref{4.1}, ${}_H \FF G {}_1
\cong \Dia_H(1_1) \oplus 2 \, \Dia_H(1_2)$.
Also, ${}_H \FF H {}_{C'_2} \cong
\Dia_H((C'_2)_1) \oplus \Dia_H
(C'_2 {\ua} 1_2)$ and ${}_H \FF H {}_R
\cong \Dia_H(R_1)$ for $R \in \SS -
\{ 1, C'_2 \}$. By Lemma \ref{2.6} and
Mackey decomposition,
${}_{H \times 1} \Res {}_{H \times C'_2}
(\Dia_H(C'_2 {\ua} 1_2)) \cong
2 \, \Dia_H(1_2)$. So ${}_{H \times 1} \Res
{}_{H \times C'_2} (\Dia_H((C'_2)_1)
\cong \Dia_H(1_1)$. Therefore, the
multiposet $\PP_H$ has the following
Hasse diagram, where single and double
lines indicate relative multiplicities $1$
and $2$, respectively.

\smallskip
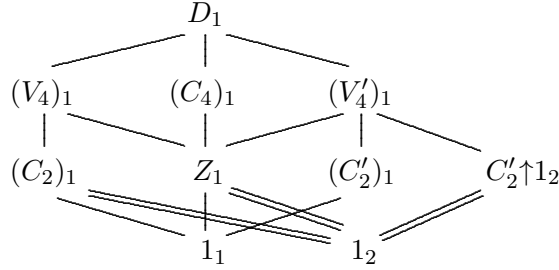
\hspace{0.9in}
\begin{picture}(300,106)

\put(102,7){$1_1$}
\put(159,7){$1_2$}
\put(30,37){$(C_2)_1$}
\put(99,37){$Z_1$}
\put(150,37){$(C'_2)_1$}
\put(210,37){$C'_2 {\ua} 1_2$}
\put(30,67){$(V_4)_1$}
\put(90,67){$(C_4)_1$}
\put(150,67){$(V'_4)_1$}
\put(97,97){$D_1$}

\curve(96,92,46,80)
\curve(104,92,104,80)
\curve(112,92,158,80)

\curve(43,62,43,50)
\curve(52,62,96,50)
\curve(104,62,104,50)
\curve(156,62,112,50)
\curve(163,62,163,50)
\curve(169,62,208,48)

\curve(48,30,96,16)
\curve(104,32,104,18)
\curve(155,30,112,16)

\curve(60,33,151,15)
\curve(114,34,156,20)
\curve(208,33,170,16)

\curve(60,31,151,13)
\curve(113,32,155,18)
\curve(209,31,171,14)

\end{picture}

Still putting $p = 2$, now put $G = S_5$.
Let $b$ be the
principal block of $\FF G$. We embed $H$
in $G$ and take $D$ and $\SS$ to be the
same as before. Let $\PP_b$ be the
multiposet of substantive generalized
pieces of $\FF G b$ having the form
$P {\ua} Q_\phi$ where $P, Q \in \SS$. The
elements of $\PP_b$ comprise a set of
representatives for the $D$-conjugacy
classes of substantive generalized pieces
of $\FF G b$ having the form
$P {\ua} Q_\phi$ where $P \leq D$. So,
as before, to specify the poset
$\cP^|(\FF G b)$, it will be enough to
display a Hasse diagram for $\PP_b$,
labelled with the relative multiplicities.

Replacing $H$ with $G$, we retain the
above convention for expressing pieces
and generalized pieces concisely. We can
choose the enumeration $\Irr(\FF G) =
\{ \theta_1^G, \theta_2^G, \theta_3^G \}$
such that $\theta_3^G$ is not in $b$. The
irreducible $\FF G$-characters $\theta_2^G$
and $\theta_3^G$ both have degree $4$.
The local pointed group corresponding
to the piece $1_3$ is not on $\FF G b$,
so $1_3$ is not a piece of $\FF G b$.
Writing $C = C_G(C'_2)$,
then $C \cong C_2 {\times} S_3$. We
have $\Irr(\FF C b_C) = \{ \theta_1^C,
\theta_2^C \}$. As before, the local
pointed group corresponding to
$(C'_2)_2$ is not on the principal
block algebra of $\FF C$, so
$(C'_2)_2$ is not a piece of $\FF G b$.
Given $R \in \SS - \{ 1, C'_2 \}$, then
$C_G(R)$ is a $p$-group. So the
pieces in $\PP_b$ are $1_2$ and
$P_1$ with $P \in \SS$.

We claim that the other elements of
$\PP_b$ are $P {\ua} 1_2$, where
$P \in \{ C_2, C'_2, Z, C_4, V'_4 \}$. To
prove the claim, we shall consider, for
each $P \in \SS$, the restriction to $P$
of the simple $\FF G$-module
$V = V(\theta_2^G)$. Let $\Omega$ be
the $G$-set of Sylow $5$-subgroups of
$G$. Observe that, as a $D$-set by
restriction,
$$\Omega \cong D/C_4 \sqcup D/C_2 \; .$$
Enumerate $\Omega = \{ L_1, ..., L_6 \}$. Let
$A = \{ {\sum}_i \lambda_i L_i \in \FF \Omega :
\sum_i \lambda_i = 0 \}$ and
$B = \FF \sum_i L_i$. The $\FF G$-modules
$\FF \Omega / A$ and $B$ are trivial and the
Brauer character of $\FF \Omega$ is easily
shown to be $2 \theta_1^G + \theta_2^G$.
Therefore, $A/B \cong V$. We have
$$m(P {\ua} 1_2) = \dim_\FF(P^+ . V)
  = \dim_\FF((P^+ . A + B)/B)$$
where $P^+$ denotes the sum of the
elements of $P$.

If $P = C'_2$, then we can choose the
enumeration $L_1$, $...$ such that the
$P$-orbits are $\{ L_1, L_2 \}$,
$\{ L_3, L_4 \}$, $\{ L_5, L_6 \}$,
whereupon
$$P^+ . A = \sspan_\FF \{
  L_1 + L_2 + L_5 + L_6,
  L_3 + L_4 + L_5 + L_6 \}$$
while $B \cap P^+ . A = \{ 0 \}$, hence
$m(P {\ua} 1_2) = 2$. If $P \in \{ C_2, W \}$,
then we can choose the enumeration such
that the $P$-orbits are $\{ L_1, L_2 \}$,
$\{ L_3, L_4 \}$, $\{ L_5 \}$, $\{ L_6 \}$,
whereupon
$$P^+ . A = \sspan_\FF \{
  L_1 + L_2, L_3 + L_4 \}$$
and again $B \cap P^+ . A = \{ 0 \}$, hence
$m(P {\ua} 1_2) = 2$. If $P \in
\{ C_4, V'_4 \}$, then we can choose the
enumeration such that $\{ L_1, L_2,
L_3, L_4 \}$ is a $P$-orbit, whereupon
$$P^+ . A = \sspan_\FF \{ L_1 + L_2 +
  L_3 + L_4 \}$$
and yet again $B \cap P^+ . A = \{ 0 \}$,
hence $m(P {\ua} 1_2) = 1$. Finally,
when $P \in \{ V_4, D \}$, the whole of
$\FF \Omega$ is annihilated by $P^+$,
hence $m(P {\ua} 1_2) = 0$. The claim
is established and moreover, the
absolute multiplicities of the pieces in
$\PP_b$ having been clear already, we
have now determined the absolute
multiplicities of all the elements of
$\PP_b$.

By Corollary \ref{4.1},
$${}_G (\FF G b) {}_P = \left\{
  \begin{array}{ll}
  \Dia_G(1_1) \oplus 4 \, \Dia_G(1_2)
  & \mbox{\rm{if $P = 1$,}} \\
  \Dia_G(P_1) \oplus 2 \, \Dia(P {\ua} 1_2)
  & \mbox{\rm{if $|P| = 2$,}} \\
  \Dia_G(P_1) \oplus \Dia(P {\ua} 1_2)
  & \mbox{\rm{if $P \in \{ C_4, V'_4 \}$,}} \\
  \Dia_G(P_1)
  & \mbox{\rm{if $P \in \{ V_4, D \}$.}}
\end{array} \right. $$
Given any element of $\PP_b$ having
the form $P {\ua} 1_2$, Mackey
decomposition yields
$${}_{G \times 1} \Res {}_{G \times P}
  (\Dia_G(P {\ua} 1_2)) \cong |P| \,
  \Dia_G(1_2) \; .$$
It follows that, for any element of $\PP_b$
having the form $P_1$, we have
$${}_{G \times 1} \Res {}_{G \times P}
  (\Dia_G(P_1)) = \left\{ \begin{array}{ll}
  \Dia_G(1_1) \oplus 4 \, \Dia_G(1_2)
  & \mbox{\rm{if $P \in \{ V_4, D \}$,}} \\
  \Dia_G(1_1) & \mbox{\rm{otherwise.}}
  \end{array} \right.$$
Therefore, the multiposet $\PP_b$ has
the following Hasse diagram, the single
and double lines again indicating
relative multiplicities $1$ and $2$,
respectively.

\smallskip
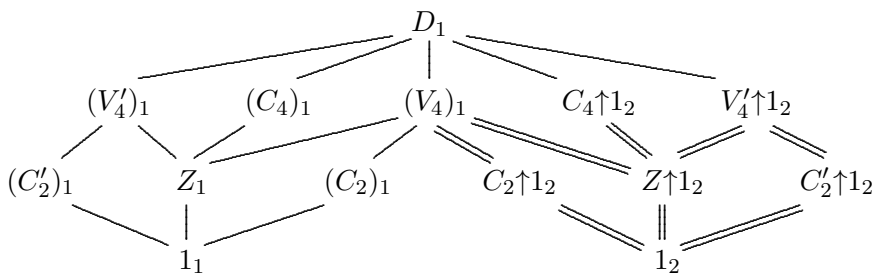
\hspace{0.0in}
\begin{picture}(350,106)

\put(85,7){$1_1$}
\put(265,7){$1_2$}
\put(20,37){$(C'_2)_1$}
\put(84,37){$Z_1$}
\put(140,37){$(C_2)_1$}
\put(200,37){$C_2 {\ua} 1_2$}
\put(260,37){$Z {\ua} 1_2$}
\put(320,37){$C'_2 {\ua} 1_2$}
\put(50,67){$(V'_4)_1$}
\put(110,67){$(C_4)_1$}
\put(170,67){$(V_4)_1$}
\put(230,67){$C_4 {\ua} 1_2$}
\put(290,67){$V'_4 {\ua} 1_2$}
\put(173,97){$D_1$}


\curve(81,16,45,31)
\curve(88,18,88,31)
\curve(95,16,140,31)

\curve(261,16,230,31)
\curve(267,18,267,31)
\curve(275,16,320,31)

\curve(260,14,229,29)
\curve(269,18,269,31)
\curve(276,14,321,29)


\curve(41,49,56,61)
\curve(84,49,70,61)
\curve(91,49,110,61)
\curve(97,47,168,64)
\curve(159,49,175,61)

\curve(204,49,183,61)
\curve(256,45,197,64)
\curve(262,49,247,61)
\curve(275,49,299,59)
\curve(330,51,310,61)

\curve(203,47,182,59)
\curve(255,43,196,62)
\curve(264,50,249,62)
\curve(274,51,298,61)
\curve(329,49,309,59)


\curve(70,79,166,95)
\curve(130,79,173,94)
\curve(180,79,180,93)
\curve(228,79,189,93)
\curve(288,79,194,95)

\end{picture}

Comparing the two Hasse diagrams
that we have produced, it may be of
interest to note that, confining attention
to the pieces, we see that the poset of
local pointed groups on a source algebra
of $\FF H$ is not isomorphic to the
poset of local pointed groups on a
source algebra of $\FF G b$.

\end{document}